
\documentstyle{amsppt}
\magnification=1200
\voffset=2cm
\def\nologo{\expandafter\redefine\csname logo\string @\endcsname{}}



\def\h5{\hskip 5pt}
\def\on{\operatorname}
\def\be{\beta}

\def\et{\eta}
\def\ze{\zeta}
\def\si{\sigma}
\def\ta{\tau}
\def\om{\omega}
\def\fracd#1#2{{\displaystyle{#1}\over\displaystyle{#2}}}
\def\Si{\Sigma}

\def\@{\char "40}
\topmatter
\title On triple coverings of irrational curves
\endtitle
\rightheadtext{On triple coverings of irrational curves}
\author Takao Kato, Changho Keem and Akira Ohbuchi
\endauthor
\leftheadtext{Takao Kato, Changho Keem and Akira Ohbuchi}
\address
Department of Mathematics,
Faculty of Science,
Yamaguchi University,
Yamaguchi, 753 Japan
\endaddress
\email
kato{\@}ccyi.ccy.yamaguchi-u.ac.jp
\endemail
\address
Department of Mathematics,
Seoul National University,
Seoul 151-742 Korea
\endaddress
\email
ckeem{\@}math.snu.ac.kr
\endemail
\address
Department of Mathematics,
Faculty of Integrated Arts abd Sciences,
Tokushima University,
Tokushima, 770 Japan
\endaddress
\email
ohbuchi{\@}ias.tokushima-u.ac.jp
\endemail
\thanks
This work was done under JSPS-KOSEF joint research program 1994.
The first named author was partially supported by A Grant-in-Aid for
Scientific Research, International Scientific Research Program, Joint
Research, The Ministry of Education, Science and Culture, Japan. While
this paper was prepared for publication, the second named author was
enjoying the hospitality of SFB 170 (G\"ottingen)
\endthanks
\abstract
Given a triple covering $X$ of genus $g$ of a general (in the sense of
Brill-Noether) curve $C$ of genus $h$, we show the existence of
base-point-free pencils of degree $d$ which are not composed with the
triple covering for any   $d\ge g-[{3h+1\over 2}]-1$ by utilizing some
enumerative  methods and computations.  We also discuss about the
sharpness of our main result and the so-called Castelnuovo-Severi bound by
exhibiting some examples. \endabstract  \endtopmatter
\document
\nologo
\specialhead 0. Introduction
\endspecialhead

In this paper, we investigate the problem of the existence of
base-point-free pencils of relatively low degree on a triple covering
$X$ of genus $g$ of a general curve $C$ of genus $h>0$. Such a problem is
classical and the picture is rather well known for the degree range
close to the genus $g$.  On the other hand, by a simple application of
the Castelnuovo-Severi inequality one can easily see that there does not
exist a base-point-free pencil of degree less than or equal to
${g-3h\over 2}$ other than the pull-backs from the base  curve $C$;
while for the degree beyond this range not many things have been known
about the existence of such a pencil which is not composed with the
given triple covering.  The main result of this paper is the following.

\proclaim{Theorem A} Let $X$ be a smooth algebraic curve of genus $g$,
over an algebraic closed field $k$ with $char(k)=0$, which admits a
three sheeted covering onto a general curve $C$ of genus $h\ge 1$,
$g\ge (2[{3h+1\over 2}]+1)([{3h+1\over 2}]+1)$. Then there exists a
base-point-free pencil of any degree $d\ge g-[{3h+1\over 2}]-1$ which is
not composed with the given  triple covering.
\endproclaim

This paper is organized as follows. In \S 1, we prove Theorem A by using
the enumerative methods and computations in $H^*(X_\alpha,\Bbb Q)$ of
various sub-loci of the symmetric product $X_\alpha$ of the given triple
covering $X$, while we defer proving the key lemma which are necessary to
prove Theorem A.  Specifically we compare the fundamental class of
$X^1_\alpha :=\{ D\in X_\alpha :\operatorname{dim}|D|\ge 1\}$ with the
class of all irreducible components of $X^1_\alpha$ whose general
elements correspond to pencils on $X$ with base points. This argument
works because the latter components are all induced from the base curve
$C$ and $X^1_\alpha$ has the expected dimension. In \S 2, we proceed to
prove Lemma 2.1. In proving Lemma 2.1, we need to utilize and
carry out some computations based on the  fundamental result of Rick
Miranda about the triple covering  of  algebraic varieties in general
[M].  In the last section we discuss  about the sharpness of our main
result and the so-called Castelnuovo-Severi bound by considering some
examples.

\specialhead
1. Existence of base-point-free pencils on triple coverings
\endspecialhead

The purpose of this section is to prove Theorem A. We begin with the
following dimension theoretic statement about the variety of special
linear systems on an algebraic curve;cf. [CKM] Corollary 3.1.2.

\proclaim{Lemma 1.1} Let X be an algebraic curve of genus $g$. Let
$n\in \Bbb N$, $g\ge (2n-3)(n-1)$ $($and if $n\le 2$ let $g\ge2n-1$$)$.
Assume that $\operatorname{dim}W^1_{n+1}(X)<1$. Then $W^1_d(X)$ is
equi-dimensional of the minimal possible dimension
$\rho (d,g,1):=2d-2-g$ for all $d$ such that $g-n<d\le g$.
\endproclaim

The following lemma is a special case of the
Castelnuovo-Severi bound which one can prove easily by using the adjunction
formula and the Hodge index theorem.

\proclaim{Lemma 1.2 \rm{(Castelnuovo-Severi bound)}}
Let $X$ be a smooth algebraic curve of genus $g$ which admits a triple
covering onto a curve $C$ of genus $h$. Let $\pi :X\rightarrow C$ be the
triple covering. For any integer $n\le {g-3h\over 2}$, a base-point-free
pencil $g^1_n$ $($possibly incomplete$)$ is composed with $\pi$, or
equivalently the morphism $X\rightarrow \Bbb P^1$ induced by the $g^1_n$
always factors through $\pi$.
\endproclaim

\demo{Proof of Theorem A}
Let $n=[{3h+1\over 2}]+2$. Then $n+1=[{3h+1\over 2}]+3\le {g-3h\over 2}$
by the numerical hypothesis on the genus $g$ of $X$. Assume that there
exists a $g^1_{n+1}$ on $X$. By Lemma 1.2, it follows  that  every
$g^1_{n+1}$ is composed with $\pi$ and
$$
W^1_{n+1}(X)=\pi^*W^1_m(C)+W_{n+1-3m}(X), \tag A.1
$$
for some $m$ with $n+1-3m\ge 0$.  On the other hand, since $C$ is general
$$
\operatorname{dim}W^1_m(C)=\operatorname {dim}\pi^*W^1_m(C)=\rho
 (m,h,1):=2m-2-h\ge 0.\tag A.2
$$
By (A.1) and (A.2), we have
$$
\operatorname{dim}W^1_{n+1}(X)<1.
$$
Accordingly by Lemma 1.1, $W^1_d(X)$ is equi-dimensional of the minimal
possible dimension $\rho (d,g,1):=2d-2-g$ for all $d$ such that
$g-n<d\le g$.  Furthermore, since
$$
\on{dim}W^1_d(X)=\rho(d,g,1)>\on{dim}(\pi^*W^1_m(C)+W_{d-3m}(X))
$$
for any $d\ge g-[{3h+1\over 2}]$ with $\rho (m,h,1)\ge 0$, we conclude
that there exists a base-point-free complete pencil of any degree
$d\ge g-[{3h+1\over2}]$ which is not composed with $\pi$. Thus
it remains to prove the theorem only for the case
$d=g-[{3h+1\over 2}]-1$. In order to avoid the unnecessary symbol [ ], we
consider the two cases according to the parity of $h$, the genus of
the base curve $C$.
\vskip 5pt
(a) {\it $h$ is even}: Set $h=2e$, then $d=g-[{3h+1\over 2}]-1=g-3e-1$.
Let $\Sigma$ be a component of $W^1_{g-3e-1}(X)$ whose general element
has a base point. Then $\Sigma = \Sigma^1_{\beta}+W_{g-3e-1-\beta}(X)$
for some $\beta$ with $4\le\beta\le g-3e-2$, where $\Sigma^1_\beta$ is a
subvariety of $W^1_\beta (X)$ whose general element is base-point-free.
Note that $\on{dim}\Sigma = \rho (g-3e-1,g,1)=g-6e-4$, hence
$\be -3e-3=\on{dim}\Sigma^1_\be\ge 0.$

Let $L\in\Sigma^1_\be$ be a general element. By the standard description
of the Zariski tangent space to the variety $W^r_d(X)$ of any curve $X$
in general, we have
$$
\hbox{\rm dim}(\operatorname{Im}\mu _0)^{\perp} =\hbox{\rm dim}
T_L(W^1_\be(X))\ge\hbox{\rm dim}\Sigma ^1_\be = \be-3e-3
$$
where
$\mu _0:H^0(X,L)\otimes H^0(X,K_X\otimes L^{-1})\rightarrow H^0(X,K_X)$ is
the usual cup product map. By the base-point-free pencil trick
([ACGH], p126), we have
$$
\eqalign{\hbox{\rm dim}(\operatorname{Im}\mu_0)^{\perp}&=g-\hbox{\rm dim}
(\operatorname{Im}\mu_0) =g-h^0(X,L)h^1(X,L)+
\hbox{\rm dim}(\operatorname{Ker}\mu_0) \cr
&=g-2(g-\be+1)+h^0(X,K_X\otimes L^{-2}) \cr
&=h^0(X,L^2)-3\ge \be -3e-3.\cr}
$$
Hence $h^0(X,L^2)\ge\be -3e\ge 3$ which implies
$\hbox{\rm dim}W^{\be -3e-1}_{2\be}(X)\ge \be -3e-3.$  By reducing to pencils
if
necessary we have
$$
\hbox{\rm dim}W^1_{\be +3e+2}(X)=
\hbox{\rm dim}W^1_{2\be -(\be -3e-2)}(X)\ge \be
-3e-3+(\be -3e-2)=2(\be -3e)-5.
$$
Note that $\be \le g-3e-2$, hence $\be +3e+2\le g$.  We now consider the
following two cases:

\noindent (1) If $\be +3e+2=g$, then by passing to residual series
$$
\hbox{\rm dim}W^1_{\be +3e+2}(X)=\hbox{\rm dim}W_{g-2}(X)=g-2\ge
2(\be -3e)-5,
$$
hence, $12e\ge g-7$ which contradts the genus assumption.

\noindent (2) If $\be +3e+2\le g-1$, we have
$$
2(\be -3e)-5\le\hbox{\rm dim}W^1_{\be +3e+2}(X)\le (\be +3e+2)-2-1,
$$
hence, $\be \le 9e+4$ by H.Martens-Mumford's theorem; cf.
[ACGH] Ch. 4, \S 5.

\vskip 5pt
On the other hand by the assumption on $h$ and $g$ we have
$\be\le 9e+4\le{g-3h\over 2}$. Hence by the Castelnuovo-Severi
inequality, every element of $\Sigma^1_\be $ is a pull-back of a
$g^1_{\be /3}$ on $C$, i.e. $\Sigma ^1_\be =\pi ^*(\Sigma^1_{\be /3}(C))$,
where $\Sigma^1_{\be /3}(C)$ is a component of $W^1_{\be /3}(C)$, i.e.
$$
\Sigma^1_\be\subset W^1_\be(X)=\pi^*W^1_{\be\over 3}(C)
$$
and $W^1_\be(X)=\pi^*W^1_{\be\over 3}(C)$ has dimension
$\rho ({\be\over 3},h,1)=2\cdot {\be\over 3}-2-2e$ since $C$ is general.  We
also have  $$
\eqalign{g-6e-4&=\on{dim}\Sigma
=\on{dim}(\Sigma^1_\be +W_{g-3e-1-\be}(X))\cr
&\le\on{dim}(W^1_\be(X)+W_{g-3e-1-\be}(X))\cr
&=(2\cdot {\be\over 3}-2-2e)+(g-3e-1-\be)\cr
&\le\on{dim}W^1_{g-3e-1}(X)=g-6e-4\quad .\cr}
$$
{}From this we get ${\be\over 3}=e+1$ and every component of
$W^1_{g-3e-1}(X)$ whose general element has a base point is of the form
$\pi^*\Sigma{'}+W_{g-6e-4}(X)$, where $\Sigma{'}$ is a component of
$W^1_{e+1}(C)$.

Since $C$ is general and the Brill-Noether number $\rho (e+1,h,1)$ is
zero, $\Si{'}$ is a single point in $g^1_{e+1}(C)$ of $W^1_{e+1}(C)$ and
there are exactly
$s:=h!\left({1\over{(h-e)!}}\cdot{1\over{(h-e+1)!}}\right)$
reduced points in $W^1_{e+1}(C)$, which is Castelnuovo's count for the number
of  $g^r_d$'s when $\rho =0$; cf. [ACGH] Theorem 1.3, p212.

Let $\tilde{\Si}{'}$ be the locus in $C^1_{e+1}$ corresponding to
$g^1_{e+1}(C)=\Si{'}$. We will show in the second section (Lemma 2.1)
that $W^1_{g-3e-1}(X)$ is reduced at a general point of
$\pi^*(W^1_{e+1}(C))+W_{g-6e-4}(X)$, which in turn says that
$X^1_{g-3e-1}$ is reduced at a general point of
$\pi^*(\tilde{\Si}{'})+X_{g-6e-4}$.

We now recall some of the notations and conventions used in [ACGH],
especially in Chapter VIII. Let $u:X_d\rightarrow J(X)$ be the abelian
sum map and let $\theta$ be the class of the theta divisor in $J(X)$.
Let $u^*:H^*(J(X),\Bbb Q)\rightarrow H^*(X_d,\Bbb Q)$ be the
homomorphism induced by $u$. By abusing notation, we use the same
letter $\theta$ for the class $u^*\theta$. By fixing a point $P$ on $X$,
one has the map $\iota:X_{d-1}\rightarrow X_d$ defined by
$\iota (D)=D+P$. We denote the class of $\iota (X_{d-1})$ in $X_d$ by
$x$.

By our discussion so far, the only components of $W^1_{g-3e-1}(X)$
whose general element has a base point is of the form
$\pi^*\Sigma{'}+W_{g-6e-4}(X)$, where $\Sigma{'}$ is a component of
$W^1_{e+1}(C)$.

We denote $\sigma$ and $\tilde{\sigma}$ by the class of
$\pi^*C^1_{e+1}+X_{g-6e-4}$ in $X_{g-3e-1}$ and the class of
$\pi^*C^1_{e+1}$ in $X_{3e+3}$ respectively. Because $X^1_{g-3e-1}$ is
of pure and expected dimension $\rho (g-3e-1,g,1)+1=g-6e-3$, the class
$x^1_{g-3e-1}$ of $X^1_{g-3e-1}$ in $X_{g-3e-1}$ is well known
(cf. [ACGH] p326), namely
$$
x^1_{g-3e-1}={1\over{(3e+1)!(3e+2)!}}\left( (3e+1)!
\theta^{3e+2}-(3e+2)!x\theta^{3e+1}\right) .\tag A.3
$$

Let's also recall that for a cycle $Z$ in $X_d$, the assignments
$$
Z\mapsto A_k(Z):=\{E\in X_{d+k}:E-D\ge 0 \h5\hbox{\rm for some}\h5 D\in Z\}
$$
$$
Z\mapsto B_k(Z):=\{E\in X_{d-k}:D-E\ge 0 \h5\hbox{\rm for some}\h5 D\in Z\}
$$
induce maps
$$
\eqalign{&A_k:H^{2m}(X_d,\Bbb Q)\rightarrow H^{2m}(X_{d+k},\Bbb Q)\cr
&B_k:H^{2m}(X_d,\Bbb Q)\rightarrow H^{2m-2k}(X_{d-k},\Bbb Q)\cr}
$$
and the so called push-pull formulas for symmetric products hold
(cf. [ACGH], p367-369).

First note that $X^1_{g-3e-1}$ and $\Si{'}$ have dimension $g-6e-3$ in
$X_{g-3e-1}$. On the other hand, by the push-pull formulas

$$
B_{g-6e-4}(x^{g-6e-3})=(g-6e-3)x
$$
and
$$
\eqalign{(\sigma\cdot x^{g-6e-3})_{X_{g-3e-1}}
&=(A_{g-6e-4}(\tilde{\sigma})\cdot x^{g-6e-3})_{X_{g-3e-1}}\cr
&=(\tilde{\sigma}\cdot
B_{g-6e-4}(x^{g-6e-3}))_{X_{3e+3}}=(g-6e-3)(\tilde{\sigma}\cdot
x)_{X_{3e+3}}\cr
&=(g-6e-3)s\cr}
$$

\noindent
by noting the fact that $(\tilde{\sigma}\cdot x)_{X_{3e+3}}=s:=
h!\left({1\over{(h-e)!}}\cdot{1\over{(h-e+1)!}}\right)$ .

By (A.3)
$$
\eqalign{\quad\quad&(x^1_{g-3e-1}\cdot x^{g-6e-3})_{X_{g-3e-1}}\cr
&=\left({1\over{(3e+1)!(3e+2)!}}{(3e+1)!
\theta^{3e+2}-(3e+2)!x\theta^{3e+1}}\right)\cdot x^{g-6e-3}\cr
&={1\over{(3e+2)!}}x^{g-6e-3}\cdot
\theta^{3e+2}-{1\over{(3e+1)!}}x^{g-6e-2}\cdot\theta^{3e+1}\cr
&={1\over{(3e+2)!}}\cdot{g!\over{(g-3e-2)!}}-{1\over{(3e+1)!}}
\cdot{g!\over{(g-3e-1)!}}\cr}
$$
where the last equality comes from the fact that
$(x^{d-\alpha}\cdot\theta^\alpha)_{X_d}=g!/(g-\alpha)!$, which is a
consequence of Poincar\'{e}'s formula; cf. [ACGH], p328.

Finally an easy calculation yields
$$
(x^1_{g-3e-1}\cdot x^{g-6e-3})_{X_{g-3e-1}}>
(\sigma\cdot x^{g-6e-3})_{X_{g-3e-1}}
$$
and this shows that there exists a component in $X^1_{g-3e-1}$ other
than those of the form  $\pi^*(C^1_{e+1})+X_{g-6e-4}$ which in turn
proves that there exists a divisor of degree $g-3e-1$ on X which moves
in a base-point-free pencil and whose complete linear series is not
composed with the given triple covering.
\vskip 5pt
We now proceed to handle the case when $h=2e+1$. Even though the method
is almost the same, we will present somewhat detailed argument and the
computations for the convenience of the reader.
\vskip 5pt
(b) {\it $h$ is odd}: Set $h=2e+1$, then $d=g-[{3h+1\over 2}]-1=g-3e-3$.
Let $\Sigma$ be a component of $W^1_{g-3e-3}(X)$ whose general element
has a base point. Then $\Sigma = \Sigma^1_{\beta}+W_{g-3e-3-\beta}(X)$
for some $\beta$ with $4\le\beta\le g-3e-4$, where $\Sigma^1_\beta$ is a
subvariety of $W^1_\beta (X)$ whose general element is base-point-free.
Note that $\on{dim}\Sigma = \rho (g-3e-3,g,1)=g-6e-8$, hence
$\be -3e-5= \on{dim}\Sigma^1_\be\ge 0.$

Let $L\in\Sigma^1_\be$ be a general element. Again by the  description
of the Zariski tangent space to the variety of special linear systems,
we have
$$
\hbox{dim}(\operatorname{Im}\mu _0)^{\perp} =\hbox{dim}T_L(W^1_\be(X))\ge
\hbox{dim}\Sigma ^1_\be\ge \be-3e-5 .
$$
By the base-point-free pencil trick, we have
$$
\eqalign{\hbox{dim}(\operatorname{Im}\mu_0)^{\perp}&=g-\hbox{dim}
(\operatorname{Im}\mu_0) =g-h^0(X,L)h^1(X,L)+
\hbox{dim}(\operatorname{Ker}\mu_0) \cr
&=g-2(g-\be+1)+h^0(X,K_X\otimes L^{-2})\cr
&=h^0(X,L^2)-3\ge \be -3e-5.\cr}
$$
Hence $h^0(X,L^2)\ge\be -3e-2\ge 3$ which implies
$W^{\be -3e-3}_{2\be}(X)\ge \be -3e-5$.  By reducing to pencils if
necessary we have
$$
\hbox{\rm dim}W^1_{\be +3e+4}(X)=\hbox{\rm dim}W^1_{2\be -(\be -3e-4)}(X)
\ge \be -3e-5+(\be -3e-4)=2(\be -3e)-9.
$$
Note that $\be \le g-3e-4$ hence $\be +3e+4\le g$.  We consider the
following two cases:

\noindent (1) If $\be +3e+4=g$, then by passing to residual series
$$
\hbox{\rm dim}W^1_{\be +3e+4}(X)=\hbox{\rm dim}W_{g-2}(X)=g-2\ge
2(\be -3e)-9,
$$
hence $12e\ge g-15$, contradictory to the genus assumption.

\noindent (2) If $\be +3e+4\le g-1$, we have
$$
2(\be -3e)-9\le
\hbox{\rm dim}W^1_{\be +3e+4}(X)\le (\be +3e+4)-2-1,
$$
hence, $\be \le 9e+10$ by H.Martens-Mumford's theorem.
\vskip 5pt
On the other hand by the assumption on $h$ and $g$ we have
$\be \le 9e+10\le {g-3h\over 2}$. Hence by the Castelnuovo-Severi
inequality, every element of $\Sigma^1_\be$ is a pull-back of a
$g^1_{\be /3}$ on $C$, i.e. $\Sigma^1_\be =\pi^*(\Sigma^1_{\be /3}(C))$,
where $\Sigma^1_{\be /3}(C)$ is a component of $W^1_{\be /3}(C)$,
$$
\Sigma^1_\be\subset W^1_\be(X)=\pi^*W^1_{\be\over 3}(C)
$$
and $W^1_\be(X)=\pi^*W^1_{\be\over 3}(C)$ has dimension
$\rho ({\be\over 3},h,1)=2\cdot {\be\over 3}-2-(2e+1)$ since $C$ is general.
We also have
$$
\eqalign{g-6e-8&=\on{dim}\Sigma =\on{dim}(\Sigma^1_\be +
W_{g-3e-3-\be}(X))\cr
&\le\on{dim}(W^1_\be(X)+W_{g-3e-3-\be}(X))\cr
&=(2\cdot {\be\over 3}-2-2e-1)+(g-3e-3-\be)\cr
&\le\on{dim}W^1_{g-3e-3}(X)=g-6e-8\cr}.
$$
{}From this we get ${\be\over 3}=e+2$ and every component of
$W^1_{g-3e-3}(X)$ whose general element has a base point is of the form
$\pi^*\Sigma{'}+W_{g-6e-9}(X)$, where $\Sigma{'}$ is a component of
$W^1_{e+2}(C)$.

Since $C$ is general and the Brill-Noether number $\rho (e+2,h,1)=1>0$,
$W^1_{e+2}(C)$ is irreducible of dimension $1$ which follows from
Gieseker's theorem by Fulton and Lazarsfeld's connectedness theorem
together with  Brill-Noether theorem by Griffiths-Harris; cf. [ACGH] p214.
Hence we have
$$
\Si{'}=W^1_{e+2}(C)\quad \hbox{\rm and}\quad \Sigma
=\pi^*W^1_{e+2}(C)+W_{g-6e-9}(X).
$$
By Lemma 2.1,
$W^1_{g-3e-3}(X)$ is reduced at a general point of $\Si$, which in turn
says that $X^1_{g-3e-3}$ is reduced at a general point of
$\pi^*C^1_{e+2}+X_{g-6e-9}$.

We denote $\sigma$ and $\tilde{\sigma}$ by the class of
$\pi^*C^1_{e+2}+X_{g-6e-9}$ in $X_{g-3e-3}$ and the class of
$\pi^*C^1_{e+2}$ in $X_{3e+6}$ respectively.  Because $X^1_{g-3e-3}$ is
of pure and expected dimension $\rho (g-3e-3,g,1)+1=g-6e-7$, the
fundamental class $x^1_{g-3e-3}$ of $X^1_{g-3e-3}$ in $X_{g-3e-3}$ is

$$
x^1_{g-3e-3}={1\over{(3e+3)!(3e+4)!}}\left( (3e+3)!
\theta^{3e+4}-(3e+4)!x\theta^{3e+3}\right) .\tag A.4
$$

First note that $X^1_{g-3e-3}$ and $\Si$ have dimension $g-6e-7$ in
$X_{g-3e-3}$. As in the even genus case, we have

$$
B_{g-6e-9}(x^{g-6e-7})={g-6e-7\choose g-6e-9}x^2
$$
and
$$
\eqalign{(\sigma\cdot x^{g-6e-7})_{X_{g-3e-3}}
&=(A_{g-6e-9}(\tilde{\sigma})\cdot x^{g-6e-7})_{X_{g-3e-3}}\cr
&=(\tilde{\sigma}\cdot
B_{g-6e-9}(x^{g-6e-7}))_{X_{3e+6}}={g-6e-7\choose g-6e-9}(\tilde{\sigma}
\cdot x^2)_{X_{3e+6}}\cr
&={g-6e-7\choose g-6e-9}(c^1_{e+2}\cdot x^2)_{C_{e+2}}\cr
&={g-6e-7\choose g-6e-9}
\left({(2e+1)!\over e!(e+1)!}-{(2e+1)!\over (e-1)!(e+2)!}\right) ,\cr}
$$
where the last equality comes from Poincar\'{e}'s formula, the fact that
$C$ is general and the class  $c^1_{e+2}$ of $C^1_{e+2}$ in $C_{e+2}$ is
$$
c^1_{e+2}={\theta^e\over e!}-{x\theta^{e-1}\over (e-1)!}.
$$

On the other hand, by (A.4) and Poincar\'{e}'s formula
$$
\eqalign{\quad\quad&(x^1_{g-3e-3}\cdot x^{g-6e-7})_{X_{g-3e-3}}\cr
&={1\over{(3e+3)!(3e+4)!}}\left( (3e+3)!
\theta^{3e+4}-(3e+4)!x\theta^{3e+3}\right)
\cdot x^{g-6e-7}\cr
&={g!(g-6e-7)\over (3e+4)!(g-3e-3)!}.\cr}
$$

Finally an easy calculation yields
$$
(x^1_{g-3e-3}\cdot x^{g-6e-7})_{X_{g-3e-3}}>
(\sigma\cdot x^{g-6e-7})_{X_{g-3e-3}}
$$
which shows that there exists a component in $X^1_{g-3e-3}$ other than
those of the form $\pi^*(C^1_{e+2})+X_{g-6e-9}$ which in turn proves
that there exists a divisor of degree $g-3e-3$ on X which moves in a
base-point-free pencil and whose complete linear series is not composed
with the given triple covering. $~\square$
\enddemo
\vskip 5pt

\specialhead
2. The proof of Lemma 2.1
\endspecialhead

In this section, our aim is to prove the following key lemma.
\vskip 5pt

\proclaim{Lemma 2.1} Assume that $C$ is a general curve of genus $h$.

\noindent\hbox{\rm (i)} Assume  $h=2e$ and $g>6h+4$. Then a general element
${\Cal
L}$  of  $\pi^*(W_{e+1}^1(C))+W_{g-6e-4}(X)$ is a reduced point of
$W^1_{g-3e-1}(X)$.

\noindent\hbox{\rm (ii)} Assume $h=2e+1$ and $g>6h+7$. Then a general element
${\Cal L}$
of  $\pi^*(W_{e+2}^1(C))+W_{g-6e-9}(X)$ is a reduced point of
$W^1_{g-3e-3}(X)$.
\endproclaim

Assume that
$\pi :X\rightarrow C$ is a finite morphism of degree $n$.
Because ${\Cal O}_C\rightarrow \pi_*{\Cal O}_X$ is an algebra
homomorphism, the trace map gives a splitting.
Hence we have $\pi_*{\Cal O}_X={\Cal O}_C\oplus {\Cal E}$
where ${\Cal E}$ 
is defined by
$$
{\Cal E}=\{ ~a\in \pi_*({\Cal O}_X)~|~a^n+c_1a^{n-2}
+\cdots +c_{n-2}a+c_{n-1}=0~\hbox{\rm and}~c_1,\cdots ,
c_{n-1}\in {\Cal O}_C~\}.
$$
This ${\Cal E}$ is clearly a locally free sheaf of rank $n-1$.
%

\remark{Remark 2.2} For the rest of this section, we assume, again, that
$n=3$, i.e. a triple  covering. Therefore $\pi$ is a finite separable
morphism of degree 3 and  ${\Cal E}$ is a locally free ${\Cal O}_C$-module
of rank 2.  As  $\pi_*({\Cal O}_X)$ is a commutative ${\Cal O}_C$-algebra,
we have a map  defined by the multiplication,
$$
\pi_*({\Cal O}_X)\otimes_{{\Cal O}_C}
\pi_*({\Cal O}_X)\rightarrow \pi_*({\Cal O}_X).
$$
Hence the multiplication induces ${\Cal O}_C$-module homomorphisms
$$
\phi_1:S^2{\Cal E}\rightarrow {\Cal O}_C
\hskip 12pt \hbox{\rm and} \hskip 12pt
\phi_2:S^2{\Cal E}\rightarrow {\Cal E}.
$$
Furthermore, the algebra structure on
$\pi_*({\Cal O}_X)={\Cal O}_C\oplus {\Cal E}$ is written as follows:

$$
(a,b)\cdot (a^\prime ,b^\prime )=(aa^\prime +\phi_1(bb^\prime ),
ab^\prime +a^\prime b+\phi_2(bb^\prime )), \tag 2.A
$$
where $a, a^\prime\in {\Cal O}_C$ and $b, b^\prime\in {\Cal E}$.
Locally, $\phi_2$ is written as the following form [M]:
$$
\phi_2(z^2)=az+bw,
\phi_2(zw)=-dz-aw,
\phi_2(w^2)=cz+dw
$$
where $z,w\in {\Cal E}$ is a local frame of ${\Cal E}$ and
$a,b,c,d\in{\Cal O}_C$.  Such a map $\phi_2$ is called a triple covering
homomorphism due to R. Miranda.  We denote
$
\hbox{\rm TCHom}_{{\Cal O}_C}(S^2{\Cal E},{\Cal E})
$
by the set of all triple covering homomorphisms on ${\Cal E}$
and there is a functorial homomorphism\par
$$
F:\hbox{\rm TCHom}_{{\Cal O}_C}(S^2{\Cal E},{\Cal E})\rightarrow
\hbox{\rm Hom}_{{\Cal O}_C}(S^2{\Cal E},{\Cal O}_C)
$$
such that $\phi_2\in\hbox{\rm TCHom}_{{\Cal O}_C}(S^2{\Cal E},{\Cal E})$
and $\phi_1=F(\phi_2)$ gives an ${\Cal O}_C$-algebra structure
on ${\Cal O}_C\oplus {\Cal E}$ by the relation (2.A);[M], p.1131
Proposition 3.5.
\endremark

\vskip 5pt
Also for the triple covering $\pi :X\rightarrow C$ there is a functorial
isomorphism
$$
\hbox{\rm Hom}_{{\Cal O}_X}(\pi^* \pi_* ({\Cal O}_X),{\Cal O}_X)\cong
\hbox{\rm Hom}_{{\Cal O}_C}(\pi_* ({\Cal O}_X),\pi_* ({\Cal O}_X)).
$$
We take a
$\sigma\in\hbox{\rm Hom}_{{\Cal O}_X}(\pi^* \pi_* ({\Cal O}_X),{\Cal O}_X)$
 which corresponds to the identity map$:\pi_* ({\Cal O}_X)
\rightarrow \pi_*({\Cal O}_X)$. Since $\pi_* ({\Cal O}_X)=
{\Cal O}_C \oplus {\Cal E}$, we can write $\sigma=(\sigma_1,\sigma_2)$
where $\sigma_1 :{\Cal O}_X\rightarrow {\Cal O}_X$ and
$\sigma_2 :\pi^* {\Cal E}\rightarrow {\Cal O}_X$. Because $\sigma$ is
an ${\Cal O}_X$-algebra homomorphism,  $\sigma_1$ is the identity map.
As for $\sigma_2$, we have the following result.

\proclaim{Lemma 2.3} The image of $\sigma_2$ is an invertible sheaf.
\endproclaim

\demo{Proof} Note that $\hbox{\rm Im}(\sigma_2)$ is an ideal sheaf of
${\Cal O}_C$.  Thus, if $\sigma_2$ is not a zero map then
$\hbox{\rm Im}(\sigma_2)$ is an invertible sheaf, because
$\hbox{\rm Im}(\sigma_2)$ is an ideal sheaf of the divisor
$D=\hbox{\rm Spec}({\Cal O}_C/\hbox{\rm Im}(\sigma_2))$ and $X$ is a
non-singular curve, $D$ is a Cartier divisor on $X$. Therefore
$\hbox{\rm Im}(\sigma_2)$ is an invertible sheaf if $\sigma_2$ is not a
zero map. Thus we need to show that $\sigma_2$ is not a zero map.

Let $U=\hbox{\rm Spec}(A)$ be an affine open subset and let
$\pi^{-1}(U)=\hbox{\rm Spec}(B)$; note that $\pi$ is a finite morphism
therefore is an affine morphism. Let $\varGamma (U, {\Cal E})=M$. The
morphism $\sigma_2$ on $\pi^{-1}(U)$ is written as follows:\par

$$\matrix
M\otimes_A B          &\overset{\sigma_2}\to\longrightarrow~~~&\!\!\! B \cr
\Vert\!\!\;\wr        &~~~\nearrow \tau& \cr
M\otimes_A(A\oplus M) &                                 &  \cr
\endmatrix$$
\noindent
where $\tau=\mu\oplus (\phi_1, \phi_2 )$, $\mu :M\otimes_A A\rightarrow B$
is a multiplication map, $\phi_2$
is the triple covering homomorphism corresponding to the $A$-algebra
structure of $B$ and $\phi_1$ is the homomorphism defined in Remark 2.2.
Clearly $\mu$ is a non-zero map, therefore $\sigma_2$ is not a zero map
either.~$\square$
\enddemo
\vskip 5pt

%
%
%

In general, one knows that to give a morphism of $X$ to the
projective bundle ${\Bbb P}({\Cal E})$ over $C$ is equivalent to give an
invertible sheaf ${\Cal L}$ on $X$ and a surjective homomorphism
$\pi^*{\Cal E}\rightarrow {\Cal L}$;cf. [H], p162 Proposition 7.12.  Hence
by Lemma 2.3, we have a morphism  $f:X\rightarrow {\Bbb P}({\Cal E})$ such
that the following diagram  commutes.\par

$$
\matrix
\!X~&\overset{f}\to\longrightarrow &~\,{\Bbb P}({\Cal E}) \cr
\,\pi\searrow &\circlearrowleft &\!\swarrow\varphi \cr
 &C &  \cr
\endmatrix
$$

With these preparations we now prove the following lemma.

\proclaim{Lemma 2.4} $f$ is embedding.
\endproclaim

\demo{Proof} Let 
$U=\hbox{\rm Spec}(A)$ be an affine open subset of $C$. Let
$Z=\pi^{-1}(U)=\hbox{Spec}(B)$,
$\varGamma (U, \pi_*({\Cal O}_X ))=A\oplus As_1\oplus As_2$ and
$Z_i=\{ P\in Z | s_i(P)\neq 0 \}$, $i=1,2$.
Because $\pi$ is a finite morphism, $Z$ is an affine open set and hence
$Z_i$ is also an affine open set for $i=1,2$.
Moreover $Z_i=\hbox{\rm Spec}(B_{s_i})$. Therefore
$B_{s_i}$ is generated by $1, \fracd{s_j}{s_i}$
where $j\neq i$.
By a criterion for a morphism to a
projective space to be an embedding ([H] p.151 Proposition 7.2),
$f|_{\pi^{-1}(U)}$ is embedding  for every
affine open  set $U\subset C$ and hence $f$ is
embedding.~$\square$\enddemo \vskip 5pt

We now consider the ruled surface
$\varphi :{\Bbb P}({\Cal E})\rightarrow C$. Let $C_0$ be a divisor on
${\Bbb P}({\Cal E})$ such that
$\varphi_*({\Cal O}_{{\Bbb P}({\Cal E})}(C_0))={\Cal E}$.
Let $\bar{C}$ be a minimal section of
$\varphi :{\Bbb P}({\Cal E})\rightarrow C$,
$\bar{\Cal E}=\varphi_*({\Cal O}_{{\Bbb P}({\Cal E})}
(\bar{C}))$ and $\delta=-(\bar{C}^2)$.

\proclaim{Lemma 2.5} There exist an exact sequence
$$
0\rightarrow {\Cal O}_C\rightarrow\bar{\Cal E}\rightarrow
\wedge^2 \bar{\Cal E}\rightarrow 0
$$\par\noindent
and an isomorphism
$$
{\Cal E}\cong \bar{\Cal E}\otimes_{{\Cal O}_C} {\Cal M}
$$\par\noindent
for some line bundle ${\Cal M}$ on $C$.
\endproclaim

\demo{Proof} Let $s\in \varGamma (C,\bar{\Cal E})$ be a non-zero
section. Then $s$ determines an injection
$0\rightarrow {\Cal O}_C\rightarrow \bar{\Cal E}$. Let ${\Cal P}$ be the
cokernel of $0\rightarrow {\Cal O}_C \rightarrow \bar{\Cal E}$.
${\Cal P}$ is of rank 1 because ${\Cal E}$ is of rank 2.  Therefore it
is enough to show that ${\Cal P}$ is a torsion free sheaf because $C$
is non-singular. Now we assume that ${\Cal P}$ has a torsion part. Let
${\Cal F}$ be the inverse image of the torsion part of ${\Cal P}$ by the
map $\bar{\Cal E}\rightarrow {\Cal P}$. As $\bar{\Cal E}$ is torsion
free, ${\Cal F}$ is a torsion free ${\Cal O}_C$-module of rank 1.
Therefore ${\Cal F}$ is an invertible sheaf and $\deg ({\Cal F})>0$
because $C$ is non-singular and ${\Cal O}_C\subsetneqq {\Cal F}$.
But this is contradictory to the definition of $\bar{\Cal E}$;
see [H] p.373 Notation 2.8.1. The remaining part is proved in [H],
p.372 Proposition 2.8.~$\square$\enddemo \vskip 5pt

We take an invertible sheaf ${\Cal M}$ such that
${\Cal E}\cong\bar{\Cal E}\otimes {\Cal M}$ and let $n=\deg ({\Cal M})$.
Then we have $C_0\sim \bar{C}+\pi^*({\Cal M})$ (linearly equivalent)
by Lemma 2.5, therefore $C_0\equiv \bar{C}+nF$ (where $\equiv$ means
numerical equivalence and $F$ is a fibre of
$\varphi$) and $(C_0^2)=2n-\delta =\deg (\wedge^2{\Cal
E})$. One also knows the following numerical equivalence;
$$
K_{{\Bbb P}({\Cal E})}\equiv -2C_0+((2h-2)+\deg (\wedge^2 {\Cal E}))F,
\tag 2.B
$$
\noindent see [H],
p.372 Proposition 2.8  and p.374 Corollary 2.11.
%
\vskip 5pt

By (2.B), we have the following.

\proclaim{Lemma 2.6} We have $n=\displaystyle{\frac{(C_0^2)+\delta}{2}}=
\displaystyle{\frac{\delta -g+3h-2}{2}}$, and we also have
$$
f(X)\equiv 3C_0-2\deg (\wedge^2 {\Cal E})F
\equiv 3\bar{C}+\frac{3\delta +g-3h+2}{2}F
$$\par\noindent
for numerical equivalence.
\endproclaim

\demo{Proof} By the Riemann-Hurwitz relation ([H], p.127 Exercise 5.16), we
have
$$
2g-2=3(2h-2)+\deg (\wedge^2 {\Cal E}^{\otimes -2})
$$
and hence
$$
\deg (\wedge^2 {\Cal E})=\displaystyle{\frac{3(2h-2)-(2g-2)}{2}}.
$$
As $C_0\equiv \bar{C}+nF$, we have
$n=\displaystyle{\frac{\delta -g+3h-2}{2}}$ because
$(C_0^2)=2n-\delta=\deg (\wedge^2{\Cal E})$. By Lemma 2.4,
$f:X\rightarrow {\Bbb P}({\Cal E})$ is embedding and
$\pi:X\rightarrow C$ is a finite morphism of degree 3. Therefore we have
$f(X)\sim 3C_0+\varphi^*(T)$ for some divisor $T$ on $C$ and
$2g-2=(K_{\Bbb P{(\Cal E)}}+f(X)).f(X)$ by the adjunction formula. Hence
$\deg (T)=-2\deg (\wedge^2 {\Cal E})$ by (2.B) and we have
$$
\eqalignno{f(X)&\equiv 3C_0-2\deg (\wedge^2 {\Cal E})F\cr
&\equiv 3(\bar{C}+\displaystyle{\frac{\delta -g+3h-2}{2}}F)
-2(\displaystyle{\frac{3(2h-2)-(2g-2)}{2})F}\cr
&\equiv 3\bar{C}+\displaystyle{\frac{3\delta +g-3h+2}{2}}F.
\quad~\square\cr}
$$
\enddemo

\proclaim{Proposition 2.7} There is an exact sequence
$$
0\rightarrow {\Cal M}\rightarrow {\Cal E}\rightarrow {\Cal L}
\rightarrow 0,
$$\par\noindent
where $\deg ({\Cal M})=-(\displaystyle{\frac{g-3h}{2}}+\displaystyle{
\frac{2-\delta}{2}})$ and $\deg ({\Cal L})=-(\displaystyle{\frac{g-3h}{2}}+
\displaystyle{\frac{2+\delta}{2}})$.
\endproclaim

\demo{Proof} By Lemma 2.5, we have
${\Cal E}\cong \bar{\Cal E}\otimes {\Cal M}$ hence
$\wedge^2 {\Cal E}\cong \wedge^2 \bar{\Cal E}\otimes{\Cal M}^{\otimes 2}$
and $\deg (\wedge^2 \bar{\Cal E})=-\delta$. Again by the Riemann-Hurwitz
relation,
$$
\deg (\wedge^2 {\Cal E})=\displaystyle{\frac{3(2h-2)-(2g-2)}{2}}.
$$
Therefore we have $\deg ({\Cal M})=-(\displaystyle{\frac{g-3h}{2}}+
\displaystyle{\frac{2-\delta}{2}})$.
Also by Lemma 2.5, we have the exact sequence
$$
0\rightarrow {\Cal O}_C\rightarrow\bar{\Cal E}\rightarrow
\wedge^2 \bar{\Cal E}\rightarrow 0.
$$\par\noindent
By taking ${\Cal L}=\wedge^2 \bar{\Cal E}\otimes {\Cal M}$,
we have the exact sequence
$$
0\rightarrow {\Cal M}\rightarrow {\Cal E}\rightarrow {\Cal L}
\rightarrow 0
$$\par\noindent
and $\deg ({\Cal L})=-(\displaystyle{\frac{g-3h}{2}}+\displaystyle{
\frac{2+\delta}{2}})$.$~\square$
\enddemo
\vskip 5pt





%

\proclaim{Lemma 2.8} $-h\le\delta\le\displaystyle{\frac{g-3h+2}{3}}$.
\endproclaim

\demo{Proof} If $\delta\le 0$, then $-h\le \delta$ by a theorem of
Nagata;[N] p.191 Theorem 1.   Therefore we may assume that $\delta\ge 0$.
As  $f(X)\equiv 3\bar{C}+\displaystyle{\frac{3\delta +g-3h+2}{2}F}$ by
Lemma 2.6 and  $f(X)$ is irreducible, we have
$
3\delta\le \frac{3\delta +g-3h+2}{2}
$
by [H] p.382 Proposition 2.20. Therefore we have
$
\delta\le \frac{g-3h+2}{3}. \square
$
\enddemo

\proclaim{Proposition 2.9}
$\deg ({\Cal M})\le \displaystyle{\frac{-g+3h-2}{3}}$
and $\deg ({\Cal L})\le \displaystyle{\frac{-g+4h-2}{2}}$.
\endproclaim

\demo{Proof} By Proposition 2.7, we have the exact sequence
$$
0\rightarrow {\Cal M}\rightarrow {\Cal E}\rightarrow{\Cal L}\rightarrow 0
$$\par\noindent
such that $\deg ({\Cal M})=-(\displaystyle{\frac{g-3h}{2}}+\displaystyle{
\frac{2-\delta}{2}})$ and $\deg ({\Cal L})=-(\displaystyle{\frac{g-3h}{2}}+
\displaystyle{\frac{2+\delta}{2}})$. Therefore we have
$$
\eqalignno{\deg ({\Cal M})&=-(\frac{g-3h}{2}+\frac{2-\delta}{2})
= -\frac{g-3h+2}{2}+\frac{\delta}{2}\cr
&\le -\frac{g-3h+2}{2}+\frac{g-3h+2}{6}\quad{\hbox{\rm by Lemma 2.8}}\cr
&=\frac{-2g+6h-4}{6}=\frac{-g+3h-2}{3}\cr}
$$\par\noindent
and
$$
\eqalignno{\deg ({\Cal L})&=-(\frac{g-3h}{2}+\frac{2+\delta}{2})
= -\frac{g-3h+2}{2}-\frac{\delta}{2}\cr
&\le -\frac{g-3h+2}{2}+\frac{h}{2}\quad\hbox{\rm by Lemma 2.8}\cr
&=\frac{-g+4h-2}{2}. \quad\square\cr}
$$
\enddemo

%
%
%
%

With all these preparations we now start to prove Lemma 2.1.

\demo{Proof of Lemma 2.1} We will only provide the proof of the case
$h=2e$ because the odd genus case can be done in the same way.

Let ${\Cal N}={\Cal O}_X(\pi^*(D)+\Delta)$ be a general element of a
component of  $\pi^*(W_{e+1}^1(C))+W_{g-6e-4}(X)$, where
${\Cal O}_X(\Delta )\in W_{g-6e-4}(X)$ is general and
${\Cal O}_C(D)\in W_{e+1}^1(C)$. Since $C$ is general,
${\Cal O}_C(D)$ is a base-point-free pencil by the Brill-Noether theorem
for general curves. Hence we deduce that ${\Cal O}_X(\pi^*(D))$ is also
base-point-free and the base locus of ${\Cal N}$ is just $\Delta$ by the
general choice of $\Delta$.
Since every component of $\pi^*(W_{e+1}^1(C))+W_{g-6e-4}(X)$ has the
expected dimension, to show the reducedness at ${\Cal N}$ it is sufficient
to show the injectivity of the Brill-Noether map
$
\varGamma (X,{\Cal N})\otimes\varGamma (X,\omega_X\otimes
{\Cal N}^{\otimes -1})\rightarrow\varGamma (X,\omega_X)
$ by the description of the Zariski tangent space to the variety
$W^r_d(X)$.
Hence by the base-point-pencil trick, it will be enough to show
that $\varGamma(X,\omega_X\otimes {\Cal N}^{\otimes -2}\otimes {\Cal
O}_X(\Delta ))=0$.

Note that  $$
\varGamma(X,\omega_X\otimes {\Cal N}^{\otimes -2}\otimes
{\Cal O}_X(\Delta ))\cong\varGamma(X,\omega_X\otimes \pi^*({\Cal O}_C
(-D)^{\otimes 2})\otimes {\Cal O}_X(-\Delta )).
$$
We will show that
$$
\dim_k(\varGamma(X,\omega_X\otimes
\pi^*({\Cal O}_C(-D)^{\otimes -2})))\le g-6e-4,
$$
which will imply
$\varGamma(X,\omega_X\otimes {\Cal N}^{\otimes -2}\otimes
{\Cal O}_X(\Delta ))=0$ by the general choice of  $\Delta$.

On the other hand, since $C$ is a general curve
${\Cal O}_C(D)$ is a reduced point of $W_{e+1}^1(C)$ (cf. [ACGH] p.214)
and hence by the base-point-free pencil trick again, we have
$$
\dim_k(\varGamma(C,{\Cal O}_C(D)^{\otimes 2}))=3.
$$
Since $\pi_*({\Cal O}_X)\cong {\Cal O}_C\oplus {\Cal E}$,
we have
$$
\pi_*(\pi^*{\Cal O}_C(D)^{\otimes 2})\cong {\Cal O}_C(2D)\oplus
({\Cal E}\otimes {\Cal O}_C(2D)).
$$
Moreover we have the exact sequence
$$
0\rightarrow {\Cal M}\otimes {\Cal O}_C(2D)\rightarrow
{\Cal E}\otimes {\Cal O}_C(2D)\rightarrow {\Cal L}\otimes {\Cal O}_C(2D)
\rightarrow 0
$$
and
$$
\deg ({\Cal M})\le \displaystyle{\frac{-g+3h-2}{3}}\quad
{\hbox{\rm ~and }}\quad \deg ({\Cal L})\le \displaystyle{\frac{-g+4h-2}{2}}
$$
\medskip\noindent by Proposition 2.7 and Proposition 2.9.
Therefore we have \par\noindent
$$
\deg ({\Cal M}\otimes {\Cal O}_C(2D))\le
\displaystyle{\frac{-g+6h+4}{3}}<0
$$
and
$$
\deg ({\Cal L}\otimes
{\Cal O}_C(2D))\le \displaystyle{\frac{-g+6h+2}{2}}<0
$$\par\noindent
by the assumption on $g$ and $h$. Therefore we have
$\varGamma (C,{\Cal E}\otimes{\Cal O}_C(2D))=0$. This implies
$\dim_k(\varGamma(X,\pi^*({\Cal O}_C(2D))))=3$ and this finishes the
proof of Lemma 2.1 by the Riemann-Roch theorem.$~\square$
\enddemo

\remark{Remark 2.10}
(i) One may show the result of Lemma 2.1 (with a slightly weaker bound on
$g$) in the following way. Here we  sketch an alternative argument for
$h=2e+1$ only; the $h=2e$ case is similar. $W^1_{e+2}(C)$ is reduced at a
general point $\Cal L$, hence $h^0(C,2\Cal L)=4$. Suppose $\pi^*\Cal L$ is
a non-reduced point on $\pi^*(W^1_{e+2}(C))=W^1_{3e+6}(X)$. Then
$h^0(X,2\pi^*\Cal L)\ge 5$ and since  $\pi^*|2\Cal L|\subsetneq
|2\pi^*\Cal L|$, it follows that $|2\pi^*\Cal L|$ is not composed with
$\pi$. Thus $X$ has a base-point-free $g^1_x$ with $x\le 6e+9$ not
composed with $\pi$ by subtracting 3 generically chosen points on $X$.
And this is contradictory to Lemma 1.2 if  $g\ge 18e+21$. Hence
$W^1_{3e+6}(X)$ is reduced at a general point $\pi^*\Cal L$. Finally, one
can assert that  $\pi^*\Cal L\otimes\Cal O_X(\Delta)$ is a reduced point
of  $W^1_{g-3e-3}(X)$ where $\Cal O_X(\Delta)\in W_{g-6e-9}$ is general,
by  invoking a result of Coppens [C], Corollary 8.

(ii) On the other hand,  one may recognize the importance of the
method (i.e. utilizing the Rick Miranda's triple covering result) adopted
in the proof of Lemma 2.1 as follows:

(1) For $d\ge g-[{3h+1\over 2}]-1$, $W^1_d(X)$ has expected dimension
and hence the enumerative method works to prove the existence of a
base-point-free pencil which is not composed with the triple covering.
On the other hand, for the degree outside this range $W^1_d(X)$ does not
have the expected dimension and one should find a totally different
way to investigate the problem. This constitutes a part of our forthcoming
paper where the methods we used in Lemma 2.1 play some role in the
paper.

(2) The reducedness result we proved in Lemma 2.1 has somewhat better bound
(of the genus $g$ of the curve $X$ with respect to the genus  $h$ of the
base curve $C$) than the one we outlined in (i) above; $g\ge 12e+14$
vs. $g\ge 18e+21$ for $h=2e+1$ and $g\ge 12e+5$
vs. $g\ge 18e+6$ for $h=2e$.
\endremark

\specialhead 3. Examples
\endspecialhead

In this section, we shall give examples of cyclic triple coverings of
curves.  One example suggests that there is a possibility of an
improvement of our bounds given in Section 1.  Another assures that the
Castelnuovo-Severi bound is not sharp in some sense, i.~e. the bound of
degree of base-point-free pencils, which is not composed with the triple
covering projection, is greater than the Castelnuovo-Severi bound.  Both
examples are given by cyclic triple coverings.  We begin with the
following remark about cyclic triple coverings.

\remark{Remark 3.1} Let $C$ be a curve of genus $h$ and $X$ be a cyclic
triple covering of $C$.  Let $X$ be of genus $g$ and let $\si$ be the
automorphism of $X$ of order 3 so that $X/\langle\si\rangle$ is
equivalent to $C$.  Let $\pi$ be the projection map of $X$ to $C$.
Since the order of $\si$ is prime, every fixed point of $\si$ is a total
ramification point of $\pi$ and they are all.  By the Riemann-Hurwitz
formula, there are $g-3h+2$ fixed points, $P_1,\dots ,P_{g-3h+2}\in X$,
of $\si$.

Since the order of $\si$ is 3, the space $H=H^1(X,{\Cal O})$ can be
decomposed into the direct sum $H=H_0\oplus H_1\oplus H_2$, where
$$
H_j=\{\om\in H:\om\circ\si =\ze\om\} ,\quad (j=0,1,2),
\quad\ze=e^{2\pi i/3}.
$$
It is known that every $\om\in H_0$ is the pullback $\pi^*\om'$ of some
$\om'\in H^1(C,{\Cal O}_C)$.  So $\dim H_0=h$ and for $\om\in H_0$, the
divisor $(\om )=2P_1+\cdots +2P_{g-3h+2}+\pi^*(Q_1+\cdots +Q_{2h-2})$,
where $Q_1+\cdots +Q_{2h-2}$ is a canonical divisor on $C$ (cf. [L]).
\endremark

For a fixed point $P_i\in X$ of $\si$, one can take a local parameter
$z_i$ at $P_i$ such that $\si :z_i\mapsto\ze z_i$ or
$\si :z_i\mapsto\ze^2 z_i$.  The exponent of $\ze$ only depends on $P_i$.

For $\et\in H_1$, let the Taylor expansion of $\et$ near $P_i$ be given
by
$$
\et =(a_0+a_1z_i+a_2z_i^2+\cdots )dz_i.
$$

In case $\si :z_i\mapsto\ze z_i$,
$$
\et\circ\si =(a_0+a_1\ze z_i+a_2\ze^2z_i^2+\dots )\ze dz_i.
$$
By the definition of $H_1$, $\et\circ\si =\ze\et$, hence we have
$a_j=0$ if $j\equiv 1, 2$ $(\hbox{\rm mod}\ 3)$.
Hence, $\et$ has a zero of order 0 modulo 3 at $P_i$.

In case $\si :z_i\mapsto\ze^2 z_i$, using a similar argument as above,
one obtains that $\et$ has a zero of order 1 modulo 3 at $P_i$.

Hence, there is an integer $t$ such that for every $\et\in H_1$,
renumbering the suffixes if necessary,
$$
(\et )=P_1+\dots +P_t+\pi^*(D_1),\tag3.1.1
$$
where $D_1$ is an effective divisor on $C$ which depends only on $\et$.
Similarly, for every $\ta\in H_2$,
$$
(\ta )=P_{t+1}+\dots +P_{g-3h+2}+\pi^*(D_2),\tag3.1.2
$$
where $D_2$ is an effective divisor on $C$ which depends only on $\ta$.

Let $k_j=\deg D_j$.  Since $\et ,\ta\in H$, $\deg\et =\deg\ta =2g-2$.
Hence, $t+3k_1=g-3h+2-t+3k_2=2g-2$.  Since $0\le t\le g-3h+2$,
$$
3k_j\ge 2g-2-(g-3h+2)=g+3h-4.
$$
If $g+3h-4\ge 6h-5$, i.~e. $k_j>2h-2$, then
$h^0(C,D_j)=\dim H_j=k_j-h+1$.
We may also assume $t\ge g-3h+2-t$, i.~e. $2t\ge g-3h+2$.  If it is not
the case, we take $\ze =e^{4\pi i/3}$ instead of $e^{2\pi i/3}$

\vskip 5pt
\proclaim{Example 3.2}Let $X, C, \pi , \si ,t$ be as above.  If
$g\ge 3h-1$, then there exists a complete base-point-free pencil
$g^1_n$ which is not composed with $\pi$, for some integer
$n\le\max\{ t, g+2-t\}$.
\endproclaim
\demo{Proof}  It is enough to show the existence of a linear series
$g^1_{\max\{ t, g+2-t\}}$ which is not composed with $\pi$.

Take an $\et\in H_1$.  Let the divisor of $\et$ be
$$
(\et )=P_1+\cdots +P_t+\pi^*(Q_1+\cdots +Q_{k_1}),
$$
$Q_1,\dots ,Q_{k_1}\in C$.  By (3.1.2), for every $\ta\in H_2$, the
divisor of $\ta$ is of the form $P_{t+1}+\cdots +P_{g-3h+2}+\pi^*(D_2)$,
$\deg D_2=k_2$ and $h^0(D_2)=k_2-h+1$.

In case $t\le\fracd{g+2}{2}$, we have $2t=g+2-3h+3k_2-3k_1\le g+2$,
hence $k_2-h\le k_1$.  So we may assume that
$$
D_2=Q_1+\cdots +Q_{k_2-h}+Q'_{k_2-h+1}+\cdots +Q'_{k_2},
$$
for $Q'_{k_2-h+1},\dots ,Q'_{k_2}\in C$.  Then, $\fracd{\et}{\ta}$ is
a meromorphic function on $X$ whose divisor is
$$
\align
P_1+\cdots +P_t & +\pi^*(Q_{k_2-h+1}+\cdots +Q_{k_1})\\
& -
(P_{t+1}+\cdots +P_{g-3h+2})-\pi^*(Q'_{k_2-h+1}+\cdots +Q'_{k_2}).%
\endalign
$$
Hence, the linear series
$|P_1+\cdots +P_t+\pi^*(Q_{k_2-h+1}+\cdots +Q_{k_1})|$ is of degree
$t+3(k_1-k_2+h)=g+2-t$ and is not composed with $\pi$.

In case $t>\fracd{g+2}{2}$, we have $k_2-h>k_1$.  In this case, we take
$$
D_2=Q_1+\cdots +Q_{k_1}+Q'_{k_1+1}+\cdots +Q'_{k_2},
$$
Then, the divisor of $\fracd{\et}{\ta}$ is
$$
P_1+\cdots +P_t-
(P_{t+1}+\cdots +P_{g-3h+2})-\pi^*(Q'_{k_1+1}+\cdots +Q'_{k_2}).
$$
Then,  $|P_1+\cdots +P_t|$ is of degree $t$ and is not composed with
$\pi$.  This completes the proof.$~\square$
\enddemo

\vskip 5pt
\noindent
\proclaim{Example 3.3} Let $X$, $C$, $\pi$, $\si$, $t$ be as above.
Assume that $g\ge 3h-1$.  If $n<\fracd{g-3h+2+t}{3}$, then every $g^1_n$
is composed with $\pi$.  Furthermore, if $g\ge 7h-4$, then for every $C$
and every $t$ such that $\fracd{g-3h+2}{2}\le t\le g-3h+2$ and
$t\equiv 2g-2$ mod $3$, there exists a cyclic triple covering $X$ of $C$
on which there exists a complete base point free linear series $g^1_t$
not composed with the triple covering.
\endproclaim
\def\ord{\hbox{\rm ord}}
\demo{Proof}  Let $f$ be a meromorphic function on $X$ such that
$f\circ\si =\ze f$, $\deg (f)_{\infty}=N_1$ and the support of
$(f)_{\infty}$ does not contain any $P_j$.  Then $f(P_j)=0$.  For
$\om\in H_0$, $\om'=f\om$ is a 1-form which satisfies
$\om'\circ\si =\ze\om'$.  Hence, the order of zero of $\om'$ at $P_j$,
denoted by $\ord_{P_j}\om'$, is
$$
\ord_{P_j}\om'\equiv\left\{
\aligned 1\quad &(j=1,\dots ,t)\cr
0\quad &(j=t+1,\dots ,g-3h+2)\endaligned\right.\quad\hbox{\rm mod}\ 3.
$$
Hence,
$$
\ord_{P_j}f\equiv\left\{
\aligned 2\quad & (j=1,\dots ,t)\cr
1\quad & (j=t+1,\dots ,g-3h+2)\endaligned\right.\quad\hbox{\rm mod}\ 3.
$$
So, $N_1=\deg (f)_{\infty}=\deg (f)_0\ge g-3h+2+t$.

Applying a similar argument for a meromorphic function $f$ on $X$ such
that $f\circ\si =\ze^2 f$, $\deg (f)_{\infty}=N_2$ and the support of
$(f)_{\infty}$ does not contain any $P_j$, we have $N_2\ge 2(g-3h+2)-t$.

Let $f$ be a meromorphic function on $X$ such that $\deg (f)_{\infty}=m$
and $f\circ\si\neq f$.  Considering a linear fraction, if necessary, we
may assume neither the support of $(f)_{\infty}$ nor $(f)_0$ contains
any $P_j$.

Let $f_i=f+\ze^{2i}f\circ\si+\ze^if\circ\si^2$, $(i=0,1,2)$.  Then
$f_i\circ\si=\ze^if_i$, $(i=0,1,2)$ and $3f=f_0+f_1+f_2$.  Since
$f(P_j)\neq 0$ and $f_i(P_j)=0$ $(i=1,2, j=1,\dots ,g-3h+2)$, we have
$f_0(P_j)\neq 0$.  Further we may assume neither $f_1$ nor $f_2$ vanishes
identically.  If it is not the case, say $f_2$ vanishes identically
(denoted by $f_2\equiv 0$), one can consider $F=\fracd{1}{f}$ instead of
$f$.  Again, $F$ can be decomposed as $F=F_0+F_1+F_2$ such that
$F_i\circ\si =\ze^iF_i$ $(i=0,1,2)$.  Assume $F_1\equiv 0$, then
$$
3=3fF=(f_0+f_1)(F_0+F_2)=(f_0F_0+f_1F_2)+f_1F_0+f_0F_2.
$$
Here, $(f_0F_0+f_1F_2)\circ\si =f_0F_0+f_1F_2$,
$f_1F_0\circ\si =\ze f_1F_0$ and $f_0F_2\circ\si =\ze^2f_0F_2$.  Hence,
$f_1F_0\equiv f_0F_2\equiv 0$ and $f_1\equiv F_2\equiv 0$.
A contradiction.  In case $F_2\equiv 0$, again we have a contradiction.

By the definition of $f_i$, we have $(f_i)_{\infty}\le
(f)_{\infty}+(f\circ\si )_{\infty}+(f\circ\si^2)_{\infty}$ $(i=0,1,2)$.
Hence, $\deg (f_i)_{\infty}\le 3m$.  This implies that
$$
3m\ge\max\{g-3h+2+t, 2(g-3h+2)-t\} .
$$
As in Example 3.2, we may assume $2t\ge g-3h+2$.  Hence,
$m\ge\fracd{g-3h+2+t}{3}$.  This implies the first assertion.

Next, we shall show the second assertion.  Let $K(C)$ be
the meromorphic function field on $C$.
Since $t\equiv 2g-2$ mod 3, $t-(g-3h+2-t)=2t-g+3h-2\equiv 0$ mod 3.  Let
$\ell=\fracd{2t-g+3h-2}{3}\ge 0$.  Take $g-3h+2-t+\ell$ general points
$P_{t+1},\dots , P_{g-3h+2},Q_1,\dots ,Q_{\ell}$ on $C$.  Since
$g-3h+2-t+3\ell =t\ge\fracd{g-3h+2}{2}\ge 2h-1$, there is a $y\in K(C)$
such that
$$
(y)_{\infty}=P_{t+1}+\cdots +P_{g-3h+2}+3(Q_1+\cdots +Q_{\ell}).
$$
We may assume the support of the zero divisor of $y$ consists of
distinct $t$ points $P_1,\dots ,P_t\in C$.
Let $x$ be a triple valued function on $C$ satisfying $y=x^3$ and let
$X$ be a curve so that $K(X)=K(C)(x)$.  For $P\in C$, taking a local
parameter $z_P$ at $P$, we have an automorphism $\si$ of $X$ such that
$\si :(z_P, x(z_P))\mapsto (z_P,\ze x(z_P))$.  Then
$X/\langle\si\rangle$ is equivalent to $C$.  Thus, $X$ is a cyclic
triple covering of $C$ which is branched over $P_1,\dots ,P_{g-3h+2}$.

Since $x\in K(X)$ satisfies $x\circ\si =\ze x$, we can divide
$P_1,\dots ,P_{g-3h+2}$ into $P_1,\dots ,P_t$ and
$P_{t+1},\dots ,P_{g-3h+2}$ as desired.  It is noted that in this case
we have $g^1_t$ which is not composed with $\pi$.  It is a better
estimate in case $t<\fracd{g+2}{2}$. $~\square$
\enddemo


\Refs
\widestnumber\key{ACGH}
\ref
\key ACGH
\by E. Arbarello, M. Cornalba, P.A. Griffiths and J. Harris
\book Geometry of Algebraic Curves I
\publ Springer Verlag
\yr 1985
\endref
\ref
\key C
\by M. Coppens
\paper Some remarks on the scheme $W^r_d$
\jour Ann. di Mat. pura ed applicata (4)
\vol 157
\pages 183-197
\yr 1990
\endref
\ref
\key CKM
\by M. Coppens, C. Keem and G. Martens
\paper Primitive linear series on curves
\jour Manuscripta Math.
\vol 77
\pages 237-264
\yr 1992
\endref
\ref
\key H
\by  R.Hartshorne
\book Algebraic Geometry
\publ Springer Verlag
\yr 1974
\endref
\ref
\key L
\by J.Lewittes
\paper Automorphisms of compact Riemann surfaces
\jour Amer. J. Math.
\vol 85
\pages 734-752
\yr 1963
\endref
\ref
\key M
\by R.Miranda
\paper Triple covers in algebraic geometry
\jour Amer. J. of Math.
\vol 107
\pages 1123-1158
\yr 1985
\endref
\ref
\key N
\by M.Nagata
\paper On self intersection number of a section on a ruled surface
\jour Nagoya Math. J.
\vol 85
\pages 191-196
\yr 1970
\endref
\endRefs
\enddocument